\documentclass[12pt]{iopart}

\usepackage{iopams}
\newtheorem{dfn}{Definition}

\newtheorem{rem}{Remark}

\newtheorem{theorem}{Theorem}  
\begin{document}

\title[ Jacobi's elliptic functions from a deformed Lie algebra ]{A novel construction of Jacobi's elliptic functions from deformed Lie algebra }

\author{Arindam Chakraborty}

\address{Physics Department,  Heritage Institute Of Technology\\
	Chowbaga Road, Anandapur, Mundapara, Kolkata-700107
	Tel.: 03366270502\\
Kolkata, West Bengal,
India\\
}
\ead{arindam.chakraborty@heritageit.edu}
\vspace{10pt}
\begin{indented}
\item[]January 2025
\end{indented}

\begin{abstract}
Jacobi's elliptic functions have been constructed from a deformed $\mathfrak{so}(2, 1)$ Lie algebra. The generators of the algebra have been obtained from a bi-orthogonal system in $\mathbb{C}^2$. The deformation parameter resembles the modulus of the relevant elliptic functions.
\end{abstract}

\noindent{\it Keywords}: Non-linear differential equation, Bi-orthogonal system, Vector fields, Lie algebra, Jacobi's elliptic functions.
%\maketitle

\section{Introduction}
Representing Lie algebras in terms of vector fields has ever been a necessity in various fields of interest ranging from classical mechanics \cite{mclach19, mukunda67, sattinger86}, quantum methods \cite{adams94}, integrable systems \cite{perelomov90, kundu03} to name a few. In such type of representation the prime objective is to find a set of $n$ functions $\{f_i : i=1\cdots n\}$ of a variable $x$ so that the set $\mathcal{L}$ of vector fields $\{V_i=f_i\partial_i : i=1\cdots n\}$ satisfy the following conditions
\begin{eqnarray}
[V_j, V_k]\in\mathcal{L}\nonumber\\
{[V_j, V_k]}= {-[V_k, V_j]}\nonumber\\
{[V_j,[V_k, V_l]]+[V_k,[V_l, V_j]]+[V_l,[V_j, V_k]]}=0
\end{eqnarray} 
where in the present case we shall define the Lie bracket $[V_j, V_k]=V_jV_k-V_kV_j$.

Elliptic functions as envisaged by Abel, Weierstrass Jacobi and many others are doubly periodic functions \cite{lawden89, akhi90} of special kind and find applications in studying real physical systems starting from finite pendulum \cite{armi06}, to orbits under central forces, motion of top \cite{lawden89}, representation of solitary wave or soliton \cite{arshad19, sarwar20, yomba23}, variable coefficient Schr\"odinger equation \cite{fan23} etc. 

In 2005 Carlson investigated Jacobi elliptic functions as inverses of a generalized form of an integral \cite{carlson05}. The integral forms of three fundamental Jacobi elliptic functions of modulus $\kappa$ have been represented as

\begin{eqnarray}
sn^{-1}(u, \kappa)&=&\int_0^u[(1-x^2)(1-\kappa^2x^2)]^{-\frac{1}{2}}dx\nonumber\\
cn^{-1}(u, \kappa)&=&\int_x^1[(1-x^2)(1-\kappa^2+\kappa^2x^2)]^{-\frac{1}{2}}dx\nonumber\\
dn^{-1}(u, \kappa)&=&\int_u^1[(1-x^2)(1+\kappa^2x^2)]^{-\frac{1}{2}}dx.
\end{eqnarray} 

Our present purpose is to suggest a possible construction of Jacobi's elliptic functions \cite{lawden89}  as compatibility conditions between a triplet of non-linear coupled ordinary differential equations involving three independent variables and the generators of a deformed version of $\mathfrak{so}(2, 1)$. The said lie algebra is obtained from a bi-orthogonal set of vectors in $\mathbb{C}^2$. Our objective is to show that the range of values of the deformation parameter plays a crucial role in the present construction and eventually resembles the modulus of the elliptic functions.

\section{Auerbach bi-orthogonal system and deformed $\mathfrak{so}(2, 1)$ algebra}
\begin{dfn}
Two pairs of vectors $\{\vert\phi_j\rangle : j=1, 2\}$ and $\vert\chi_j\rangle : j=1, 2$ in $\mathbb{C}^2$ are said to form a \textbf{bi-orthogonal set} relative to the inner-product $\langle\cdot\vert\cdot\rangle$ if $\langle\phi_j\vert\chi_k\rangle=\delta_{jk}$.
\end{dfn}

Let us consider $\vert v_j\rangle=\left(\begin{array}{c}
	c^{(1)}_j \\
	c^{(2)}_j
\end{array} \right)\in \mathbb{C}^2$ and the inner-product $\langle v_j\vert v_k\rangle=(c^{(1)}_j)^\star c^{(1)}_k+(c^{(2)}_j)^\star c^{(2)}_k$. The following theorem is immediate.

\begin{theorem}
Given a pair of vectors $\{\vert v_j\rangle : j=1,2\}$ and an invertible hermitian operator $T$ relative to the relevant inner product, the set $\{\vert\phi_j\rangle= T \vert v_j\rangle: j=1, 2\}$ and $\{\vert\chi_j\rangle= (T^{-1})^\dagger\vert v_j\rangle: j=1, 2\}$ constitute bi-orthogonal system provided
 $\langle v_j\vert v_k\rangle=0$. 
\end{theorem}

Proof : Since $T=T^\dagger$, $\langle\phi_j\vert\chi_k\rangle=\langle v_j\vert T^\dagger (T^\dagger)^{-1}\vert v_k\rangle=\langle v_j\vert v_k\rangle$. Hence, follows the theorem.$\:\:\:\square$.

Now considering
\begin{equation}\label{trans}
T=\cos\frac{\vartheta}{2}\mathbf{1}_{2\time 2}+2\cos\frac{\varphi}{2}\sin\frac{\vartheta}{2}\sigma_1-2\sin\frac{\varphi}{2}\sin\frac{\vartheta}{2}\sigma_2,
\end{equation}
where, $\sigma_1=\frac{1}{2}\left(\begin{array}{cc}
	0 & 1 \\
	1 &0
\end{array} \right)$ and $\sigma_2=\frac{1}{2}\left(\begin{array}{cc}
	0 & -i \\
	i &0
\end{array} \right)$ and the orthogonal vectors $\left\{\vert v_j\rangle=2^{-\frac{1}{2}}\left(\begin{array}{c}
	1 \\
	(-1)^{j-1}
\end{array} \right) : j=1, 2\right\}$, we obtain the following bi-orthogonal vectors

$\left\{\vert\phi_j\rangle= {\frac{1}{\sqrt{2}}}\left(\begin{array}{c}
	e^{-i(\frac{3}{2}-j)\vartheta} \\
	(-1)^{j-1}e^{i(\frac{3}{2}-j)\vartheta} 
\end{array} \right),\vert\chi_j\rangle=  {\frac{1}{\sqrt{2}}}\left(\begin{array}{c}
	e^{i(\frac{3}{2}-j)\vartheta} \\
	(-1)^{j-1}e^{-i(\frac{3}{2}-j)\vartheta} 
\end{array} \right)\right\}$ that obey $\langle\phi_j\vert\chi_k\rangle=\omega\delta_{jk}$.

Defining 
\begin{equation}
\{T^{(\gamma)}_m=\frac{1}{2}\sum_{j, k=1}^{2}\frac{\alpha_{jk}^{(m)}}{\omega^{\delta_{m3}}}\vert \phi_j\rangle\langle \chi_k\vert : m=1, 2, 3\}.
\end{equation}
with $\alpha_{jk}^{(1)}=(-i)(1-\delta_{jk})$,  $\alpha_{jk}^{(2)}=(-1)^{j-1}\delta_{jk},\alpha^{(3)}_{jk}=i(-1)^{j-1}(1-\delta_{jk})$ and $\omega=\cos\vartheta$ we get the explicit forms of $\{T_j^{(\gamma)} : j=1, 2, 3\}$ like
\begin{equation}
T^{(\gamma)}_1={\frac{1}{2}}\left(\begin{array}{cc}
	-i & \gamma\\
	\gamma & i
\end{array} \right), T^{(\gamma)}_2={\frac{1}{2}}\left(\begin{array}{cc}
	-i\gamma & 1\\
	1 & i\gamma
\end{array} \right),
T^{(\gamma)}_3={\frac{1}{2}}\left(\begin{array}{cc}
	0 & -i\\
	i & 0
\end{array} \right)
\end{equation}
as the generators of deformed $\mathfrak{so}(2, 1)$. On letting $\gamma=0$ we can recover the generators of conventional $\mathfrak{so}(2, 1)$ algebra \cite{iachello06}.

The commutation relations are given by  
\begin{eqnarray}\label{commu}
[T^{(\gamma)}_j, T^{(\gamma)}_k]=(-1)^{\delta_{l1}}(1-\gamma^2\delta_{l3})\epsilon_{jkl}T^{(\gamma)}_l.
\end{eqnarray}

\section{Differential realization of deformed $\mathfrak{so}(2, 1)$ and Jacobi elliptic functions}

Considering $\{T_j^{(\gamma)}=f_j(u)\frac{d}{du}\}$ let us compare the triplet of differential equations
 \begin{eqnarray}\label{de11}
\frac{1}{K_l^2}\frac{df_l}{du}+f_jf_k=0 :\:\: j,k,l=1,2,3\:\: (j\neq k\neq l)
\end{eqnarray}
with equations-\ref{commu} we get for example $(K_3^2f_2^2-K_2^2f_3^2)=1$ which in turn gives {\footnote{A generalized version of the differential equations considered in \cite{chul98} in the contexts of Jacobi elliptic algebra and $\mathbf{Z_2}$ grading.}}
$f_2^2+f_3^2=1\:\:{\rm{for}}\:\:K_2^2=-1, K_3^2=1$
leading to $f_3=\pm\sqrt{1-f_2^2}$. Similar consideration yields $f_1=\pm\sqrt{1-\kappa^2f_2^2}$ for $K_1^2=\kappa^2=\gamma^2$. Using these results in equation-\ref{de11} (with $j=3, k=1, l=2$) one can get
\begin{equation}
\frac{df_2}{du}=\sqrt{1-\kappa^2f_2^2}\sqrt{1-f_2^2}. 
\end{equation}

which upon
integration from $0$ to $u$ yields $f_2(u; \kappa)=sn(u; \kappa)$. $\kappa$ is called the modulus of the elliptic function. For elliptic functions with real nome, $0\leq\kappa<1$  \cite{lawden89}. Similar considerations with integration done from $1$ to $u$ give $f_1(u; \kappa)=-dn(u; \kappa)$ and $f_3(u; \kappa)=-cn(u; \kappa)$ (negative signs are chosen to fit the commutation relations). Since $f_1(u; 0)=-1, f_2(u; 0)=\sin u, f_2(u; 0)=-\cos u$ \cite{lawden89} the differential realization remains to be consistent with the matrix realization.
The Casimir of the deformed algebra is found to be $\mathcal{C}_T=T_1^2-T_2^2-(1-\kappa^2)T_3^2=\kappa^2[(cn^2u-sn^2u)\frac{d^2}{du^2}-2snu\: cnu\: dnu\frac{d}{du}]$. 

\begin{rem}
Since, $\gamma=\sqrt{1-\omega^2}=\sin\vartheta$,  $0\leq\vert\gamma\vert\leq 1$. The parameter $\gamma$ can well be replaced by $\vert\gamma\vert$. \textbf{When $\vert\gamma\vert=1$, the above bi-orthogonality fails to hold and the said algebra loses its semi-simplicity. The above generators $\{T_m^{(\gamma)}:m=1,2,3\}$ can be considered as generators of a deformed $\mathfrak{so}(2, 1)$ algebra in the sense that for $\gamma=0$, we retrieve the conventional $\mathfrak{so}(2, 1)$ generators \cite{iachello06}. On the other hand, for $\gamma=0$, the said bi-orthogonal set degenerates into orthogonal set and the elliptic functions reduce to trigonometric functions. Hence one can justify the identification of $\gamma$ with $\kappa$}. Similar generators have been used elsewhere \cite{chakraborty20} by the present author in the context of non-hermitian quantum mechanics and fusion polynomial Lie algebras.
\end{rem}

\section{Conclusion}

In view of the above discussion it may be said that the integral forms of the Jacobi's elliptic functions are closely related to the concurrence of three facts : (i) condition of bi-orthogonality of the chosen pairs of vectors, (ii) condition of semi-simplicity of the deformed Lie algebra taken into consideration and (iii) identification of the deformation parameter with the modulus of the said elliptic functions. Since the above result establishes a connection between Lie algebras and 
elliptic functions, similar method may be applied to construct functions as compatibility conditions between a set of differential equations and Lie algebra. To the best of the author's knowledge this is for the first time construction of  elliptic functions has been envisaged from the Lie algebraic perspective.
\section*{ORCID iDs}
Arindam Chakraborty https://orcid.org/0000-0002-3414-3785

\end{document}